\documentclass[12pt, reqno,fleqn]{amsart}
\usepackage{amsmath, amsthm, amscd, amsfonts, amssymb, graphicx, color}
\usepackage[bookmarksnumbered, colorlinks, plainpages]{hyperref}

\newtheorem{theorem}{Theorem}[section]
\newtheorem{lemma}[theorem]{Lemma}
\newtheorem{proposition}[theorem]{Proposition}
\newtheorem{corollary}[theorem]{Corollary}
\theoremstyle{definition}
\newtheorem{definition}[theorem]{Definition}

\theoremstyle{remark}
\newtheorem{remark}[theorem]{Remark}
\numberwithin{equation}{section}

\begin{document}

\begin{center}
\textbf{Approximation properties of Bernstein singular integrals in variable
exponent Lebesgue spaces on the real axis}

\bigskip

\textbf{Ramazan Akg\"{u}n}

\bigskip
\end{center}

\begin{quotation}
Balikesir University, Faculty of Arts and Sciences, Department of
Mathematics, \c{C}a\u{g}\i \c{s} Yerle\c{s}kesi, 10145, Bal\i kesir, T\"{u}%
rkiye,\quad \quad \newline
rakgun@balikesir.edu.tr
\end{quotation}

\bigskip

\begin{quotation}
\textbf{Abstract} In generalized Lebesgue spaces $L^{p\left( \cdot \right) }$
with variable exponent $p\left( \cdot \right) $ defined on the real axis, we
obtain several inequalities of approximation by integral functions of finite
degree. Approximation properties of Bernstein singular integrals in these
spaces are obtained. Estimates of simultaneous approximation by integral
functions of finite degree in $L^{p\left( \cdot \right) }$ are proved.

\bigskip

\end{quotation}

\bigskip

\begin{quotation}
\end{quotation}

\bigskip

\begin{quotation}
Author was supported by Balikesir University Scientific Research Project
2019/61.
\end{quotation}

\section{Introduction}

In this work we consider approximation properties of Bernstein's singular
integrals for functions given in variable exponent Lebesgue spaces $%
L^{p\left( x\right) }\left( \mathbb{R}\right) $. This scale of function
spaces were studied in details in books Uribe-Fiorenza \cite{UF13}, Diening,
Harjulehto, H\"{a}st\"{o}, R\r{u}\v{z}i\v{c}ka \cite{DHHR11} and
Sharapudinov \cite{Sh12}. $L^{p\left( x\right) }\left( \mathbb{R}\right) $
has many applications in several branches of mathematics such as elasticity
theory \cite{vz86}, fluid mechanics \cite{mr00}, \cite{krr-mr-96},
differential operators \cite{mr00}, \cite{ld-mr02}, nonlinear Dirichlet
boundary value problems \cite{zok-jr91}, nonstandard growth \cite{vz86} and
variational calculus. Variable exponent works started with W. Orlicz \cite%
{worl31} and developed in many directions. For example, $L^{p\left( x\right)
}\left( \mathbb{R}\right) $ is a modular space (\cite{jm83}) and under the
condition $p^{+}:=esssup_{x\in \mathbb{R}}p\left( x\right) <\infty $, $%
L^{p\left( x\right) }\left( \mathbb{R}\right) $ becomes a particular case of
Musielak-Orlicz spaces \cite{jm83}. Starting from nineties, studies on $%
L^{p\left( x\right) }\left( \mathbb{R}\right) $ has reached a positive
momentum: \cite{zok-jr91}, \cite{sgs94}, \cite{xf-dz01}, \cite{ld04} and
many others.

In variable exponent Lebesgue spaces on $\left[ 0,2\pi \right] $ (or $\left[
0,1\right] $), some fundamental results corresponding to the approximation
of function have been obtained by Sharapudinov \cite{Sh1, Sh2, Sh3,Sh4,Sh5}.
Some results on approximation in $L^{p\left( x\right) }\left( \left[ 0,2\pi %
\right] \right) $ or other function classes can be seen e.g. in \cite%
{acis,ra11u,ahak,AK1,AK2,AK3,ayey,day,GI2,Israfil,it16,iy,Ja,Ja1,ja2,ja3,k,ssv}%
.

In this work, we aim to obtain simultaneous theorems on approximation by
entire functions of finite degree in variable exponent Lebesgue spaces on
the whole real axis $\mathbb{R}$.

The approximation by entire function of finite degree in the real axis
started by the works of Bernstein \cite{Ber,B46}, N. Wiener and R.Paley \cite%
{RN}, N.I.Ahiezer \cite{Ak}, Nikolskii \cite{Ni}. Note that an entire
function of finite exponential type is merely an entire function of order 1
and finite type that in approximation theory, these often play an important
role similar to trigonometric polynomials in the case of approximation of
periodic functions.

Note that, some results on approximation by entire integral functions of
finite degree were obtained by Ibragimov \cite{II3} and Taberski \cite%
{T81,T86} in the classical Lebesgue spaces $L_{p}\left( \mathbb{R}\right) $.

We can give some required definitions. We denote by $P$ the class of
exponents $p(x):\mathbb{R}\rightarrow \lbrack 1,\infty )$ such that $p(x)$
is a measurable function and $p(x)$ satisfy conditions%
\begin{equation}
1\leq p_{-}:=essinf_{x\in \mathbb{R}}p(x)\text{,\quad }p^{+}:=esssup_{x\in 
\mathbb{R}}p\left( x\right) <\infty .  \label{gag21}
\end{equation}

We define the $L^{p(\cdot )}:=L^{p(\cdot )}({\mathbb{R}})$ as the set of all
functions $f:\mathbb{R}\rightarrow \mathbb{C}$ such that%
\begin{equation}
I_{p(\cdot )}\left( \frac{f}{\lambda }\right) :=\int_{\mathbb{R}}\left\vert 
\frac{f(y)}{\lambda }\right\vert ^{p(y)}dy<\infty  \label{Ip}
\end{equation}%
for some $\lambda >0$. The set of of functions $L^{p(\cdot )}$, with norm 
\begin{equation*}
\Vert f\Vert _{p(\cdot )}:=\inf \left\{ \eta >0:I_{p(\cdot )}\left( \frac{f}{%
\eta }\right) <1\right\}
\end{equation*}%
is Banach space.

For $p\in P$ we define its conjugate $p^{\prime }(x):=\frac{p(x)}{p(x)-1}$
for $p(x)>1$ and $p^{\prime }(x):=\infty $ for $p(x)=1$.

For $i\in \mathbb{N}$, all constants $\mathbf{c}_{i}$ (or $\mathbf{c}$) will
be some positive numbers such that $\mathbf{c}_{i}$ will depend on main
parameters of the problem. In some cases we will use temporaryly some
generic constans $C,c>0$ for clarity (for example in statements of some
theorems).

Throughout this paper symbol $\mathfrak{A}\lesssim \mathcal{B}$ will mean
that there exists a constant C depending only on unimportant parameters in
question such that inequality $\mathfrak{A}\leq $C$\mathcal{B}$ holds.

We will use symbol $C$ for generic constants that does not depend on main
parameters and changes with placements. We will give explicit constants in
the proofs but these constants are not best constants.

\begin{definition}
\label{Def1}Let $P^{Log}$ be a subclass (\cite{DHHR11}) of $P$ such that
there exist constants $\mathbf{c}_{1},\mathbf{c}_{2}>0$, $\mathbf{c}_{3}\in 
\mathbb{R}$ with properties%
\begin{equation}
|p(x)-p(y)|\ln \left( e+1/|x-y|\right) \leq \mathbf{c}_{1}<\infty \text{%
,\quad }\forall x,y\in \mathbb{R}\text{,}  \label{gag22}
\end{equation}%
\begin{equation}
|p(x)-\mathbf{c}_{3}|\ln \left( e+|x|\right) \leq \mathbf{c}_{2}<\infty 
\text{,\quad }\forall x\in \mathbb{R}.  \label{gag23}
\end{equation}
\end{definition}

\section{Transference Result}

Let $C_{0}^{\infty }$ be class of infinitely times continuously
differentiable functions $\phi $ with compact support $spt\phi $ in $\mathbb{%
R}$. We denote For given $f\in L^{p(\cdot )}$ we can define an auxiliary
function $F_{f}$ as follows: Define%
\begin{equation}
F_{f}\left( u\right) \text{:}\text{=}\int\nolimits_{\mathbb{R}}f\left(
u+x\right) \left\vert G\left( x\right) \right\vert dx,\quad u\in \mathbb{R},
\label{efef}
\end{equation}%
where $G\in L^{p^{\prime }(\cdot )}\cap C_{0}^{\infty }$ and $\left\Vert
G\right\Vert _{p^{\prime }(\cdot )}\leq 1.$

Let $C(A)$ be the class of continuous functions defined on $A$. We set $%
\mathbf{c}_{0}$:=$\left\Vert G\right\Vert _{\infty }$.

\begin{theorem}
\label{tra}Let $p\in P^{Log}$ and $f$,$g\in L^{p(\cdot )}$. If%
\begin{equation*}
\left\Vert F_{f,G}\right\Vert _{C\left( \mathbb{R}\right) }\lesssim
\left\Vert F_{g,G}\right\Vert _{C\left( \mathbb{R}\right) },
\end{equation*}%
with an absolute positive constant, then, we have following norm inequality%
\begin{equation*}
\left\Vert f\right\Vert _{p\left( \cdot \right) }\lesssim \left\Vert
g\right\Vert _{p\left( \cdot \right) }
\end{equation*}%
with a positive constant depending only on $p$.
\end{theorem}

\section{Mollifiers and forward Steklov means in $L^{p(\cdot )}$}

\begin{definition}
Suppose that $0<\delta <\infty $ and $\tau \in \mathbb{R}$. We define family
of translated Steklov operators $\{\mathcal{S}_{\delta ,\tau }f\}$, by%
\begin{equation}
\mathcal{S}_{\delta ,\tau }f(x):=\frac{1}{\delta }\int\nolimits_{x+\tau
-\delta /2}^{x+\tau +\delta /2}f\left( t\right) dt,\quad x\in \mathbb{R}
\label{steklR}
\end{equation}%
for$\mathcal{\ }$locally integrable function $f$ defined on $\mathbb{R}$.
\end{definition}

Let $f$ and $g$ be two real-valued measurable functions on $\mathbb{R}$. We
define the convolution $f\ast g$ of $f$ and $g$ by setting $(f\ast
g)(x)=\int\nolimits_{\mathbb{R}}f(y)g(x-y)dy$ for $x\in \mathbb{R}$ for
which the integral exists in $\mathbb{R}$.

The following result on mollifiers in variable exponent Lebesgue spaces is
obtained by D. Cruz-Uribe and A. Fiorenza (see \cite{cuf}).

\begin{definition}
\label{ddd}Let $\phi \in L_{1}\left( \mathbb{R}\right) $ and $\int\nolimits_{%
\mathbb{R}}\phi \left( t\right) dt=1$. For aech $t>0$ we define $\phi
_{t}\left( x\right) =\frac{1}{t}\phi \left( \frac{x}{t}\right) $. Sequence $%
\left\{ \phi _{t}\right\} $ will be called approximate identity. A function%
\begin{equation*}
\phi ^{\symbol{126}}\left( x\right) =\sup\limits_{\left\vert y\right\vert
\geq \left\vert x\right\vert }\left\vert \phi \left( y\right) \right\vert
\end{equation*}%
will be called radial majorant of $\phi .$ If $\phi ^{\symbol{126}}\in
L_{1}\left( \mathbb{R}\right) $, then, sequence $\left\{ \phi _{t}\right\} $
will be called potential-type approximate identity.
\end{definition}

\begin{theorem}
\label{bt}(\cite{cuf}) Suppose $p\in P^{Log}$, $f\in L^{p(\cdot )}$, $\phi $
is a potential-type approximate identity. Then, for any $t>0$,%
\begin{equation*}
\left\Vert f\ast \phi _{t}\right\Vert _{p\left( \cdot \right) }\lesssim
\left\Vert f\right\Vert _{p\left( \cdot \right) }
\end{equation*}%
and%
\begin{equation*}
\underset{t\rightarrow 0}{\lim }\left\Vert f\ast \phi _{t}-f\right\Vert
_{p\left( \cdot \right) }=0
\end{equation*}%
hold with a positive constant depend on $p$.
\end{theorem}

As a corollary of Theorem \ref{tra} we have

\begin{theorem}
\label{stekRR}Suppose that $p\in P^{Log}$, $0<\delta <\infty $ and $\tau \in 
\mathbb{R}$. Then, the family of operators $\{\mathcal{S}_{\delta ,\tau }f\}$%
, defined by (\ref{steklR}), is uniformly bounded (in $\delta $ and $\tau $)
in $L^{p(\cdot )}$, namely, for any $0<\delta <\infty $ and $\tau \in 
\mathbb{R}$ norm inequality%
\begin{equation*}
\left\Vert \mathcal{S}_{\delta ,\tau }f\right\Vert _{p\left( \cdot \right)
}\lesssim \left\Vert f\right\Vert _{p\left( \cdot \right) }
\end{equation*}%
holds with a positive constant depend on $p$.
\end{theorem}

As a corollary of Theorem \ref{stekRR} we get

\begin{corollary}
\label{coroL} Let $p\in P^{Log}$, $0<\delta <\infty $, $f\in L^{p(\cdot )}.$
If $\tau =\delta /2$ then,%
\begin{equation*}
T_{\delta }f\left( x\right) :=\mathcal{S}_{\delta ,\delta /2}f\left(
x\right) =\frac{1}{\delta }\int\nolimits_{0}^{\delta }f\left( x+t\right) dt
\end{equation*}%
and%
\begin{equation*}
\left\Vert T_{\delta }f\right\Vert _{p\left( \cdot \right) }\lesssim
\left\Vert f\right\Vert _{p\left( \cdot \right) }
\end{equation*}%
holds with a positive constant depend on $p$.
\end{corollary}

\section{Modulus of smoothness and K-functional}

If $f\in L^{p\left( \cdot \right) }$ and $0\leq \delta <\infty $, then%
\begin{equation}
\Omega _{r}\left( f,\delta \right) _{p\left( \cdot \right) }=\left\Vert
\left( I-T_{\delta }\right) ^{r}f\right\Vert _{p\left( \cdot \right)
}\lesssim \left\Vert f\right\Vert _{p\left( \cdot \right) }.  \label{Bound}
\end{equation}%
Here $I$ is the identity operator. Here and in what follows $W_{r}^{p\left(
\cdot \right) }$, $r\in \mathbb{N}$, will be the class of functions $f\in
L^{p\left( \cdot \right) }$ such that $f^{\left( r-1\right) }$ is absolutely
continuous and $f^{\left( r\right) }\in L^{p\left( \cdot \right) }$.

\begin{remark}
For $p\in P^{Log}$, $f$, $g\in L^{p\left( \cdot \right) }$ and $0\leq \delta
<\infty ,$ the modulus of smoothness $\Omega _{r}\left( f,\delta \right)
_{p\left( \cdot \right) }$, has the following usual properties:

(i) $\Omega _{r}\left( f,\delta \right) _{p\left( \cdot \right) }$\ is
non-negative; non-decreasing function of $\delta \geq 0$;

(ii) $\Omega _{r}\left( f+g,\cdot \right) _{p\left( \cdot \right) }\leq
\Omega _{r}\left( f,\cdot \right) _{p\left( \cdot \right) }+\Omega
_{r}\left( g,\cdot \right) _{p\left( \cdot \right) }$;

(iii) $\lim_{\delta \rightarrow 0^{+}}\Omega _{r}\left( f,\delta \right)
_{p\left( \cdot \right) }=0$;

(iv) $\Omega _{r}\left( f,\delta \right) _{p\left( \cdot \right) }\lesssim
\delta ^{r}\left\Vert f^{\text{ }\left( r\right) }\right\Vert _{p\left(
\cdot \right) }$ for $r\in \mathbb{N}$, $f\in W_{r}^{p(\cdot )}$ and $\delta
>0$.
\end{remark}

Indeed: (ii) follows from definition. (iii) is follow from (\ref{Bound}) and
(3.4) Theorem 3.1 of \cite{AG}. (iv) follows from Lemma 3.2 of \cite{AG}.
(i) follows from Lemma \ref{bukun} given below.

\begin{definition}
Define, for $f\in L^{p\left( \cdot \right) }$, $p\in P^{Log}$, and $\delta
>0,$%
\begin{equation*}
\left( \mathfrak{R}_{\delta }f\right) \left( \cdot \right) :=\frac{2}{\delta 
}\int\nolimits_{\delta /2}^{\delta }\left( \frac{1}{h}\int\nolimits_{0}^{h}f%
\left( \cdot +t\right) dt\right) dh.
\end{equation*}
\end{definition}

\begin{remark}
Note that, for $0<\delta <\infty $, $p\in P^{Log}$ we know from Corollary %
\ref{coroL} that%
\begin{equation*}
\left\Vert \mathfrak{R}_{\delta }f\right\Vert _{p\left( \cdot \right)
}\lesssim \left\Vert f\right\Vert _{p\left( \cdot \right) }
\end{equation*}%
and, hence, $f-\mathfrak{R}_{\delta }f\in L^{p\left( \cdot \right) }$ for $%
f\in L^{p\left( \cdot \right) }.$
\end{remark}

We set $\mathfrak{R}_{\delta }^{r}f:=\left( \mathfrak{R}_{\delta }f\right)
^{r}.$

\begin{lemma}
\label{bukun} Let $0<h\leq \delta <\infty $, $p\in P^{Log}$ and $f\in
L^{p\left( \cdot \right) }$. Then%
\begin{equation}
\left\Vert \left( I-T_{h}\right) f\right\Vert _{p\left( \cdot \right)
}\lesssim \left\Vert \left( I-T_{\delta }\right) f\right\Vert _{p\left(
\cdot \right) }  \label{bukunn}
\end{equation}%
holds with a positive constant depend on $p$.
\end{lemma}

\begin{lemma}
\label{bukunA} Let $0<\delta <\infty $, $p\in P^{Log}$ and $f\in L^{p\left(
\cdot \right) }$. Then%
\begin{equation*}
\left\Vert \left( I-\mathfrak{R}_{\delta }\right) f\right\Vert _{p\left(
\cdot \right) }\lesssim \left\Vert \left( I-T_{\delta }\right) f\right\Vert
_{p\left( \cdot \right) }
\end{equation*}%
holds with a positive constant depend on $p$.
\end{lemma}

\begin{remark}
Note that, the function $\mathfrak{R}_{\delta }f$ is absolutely continuous
and differentiable a.e. (almost everywhere) on $\mathbb{R}$ (see \cite[(5.2)
of Theorem 4]{Sh3}).
\end{remark}

The following lemma is obvious from definitions.

\begin{lemma}
\label{akgunArx} Let $0<\delta <\infty $, $p\in P^{Log}$ and $f\in
W_{1}^{p\left( \cdot \right) }$. Then%
\begin{equation}
\frac{d}{dx}\mathfrak{R}_{\delta }f=\mathfrak{R}_{\delta }\frac{d}{dx}f\text{%
\quad and\quad }\frac{d}{dx}T_{\delta }f=T_{\delta }\frac{d}{dx}f
\label{trev}
\end{equation}%
a.e. on $\mathbb{R}.$
\end{lemma}

\begin{lemma}
\label{lm05}Let $0<\delta <\infty $, $p\in P^{Log}$ and $f\in L^{p\left(
\cdot \right) }$ be given. Then%
\begin{equation}
\delta \left\Vert \frac{d}{dx}\mathfrak{R}_{\delta }f\right\Vert _{p\left(
\cdot \right) }\lesssim \left\Vert \left( I-T_{\delta }\right) f\right\Vert
_{p\left( \cdot \right) }  \label{func2}
\end{equation}%
holds with a positive constant depend on $p$.
\end{lemma}

The following lemma can be proved using induction on $r$.

\begin{lemma}
\label{da}Let $0<\delta <\infty $, $r-1\in \mathbb{N}$, $p\in P^{Log}$, and $%
f\in L^{p\left( \cdot \right) }$ be given. Then%
\begin{equation*}
\frac{d^{r}}{dx^{r}}\mathfrak{R}_{\delta }^{r}f=\frac{d}{dx}\mathfrak{R}%
_{\delta }\frac{d^{r-1}}{dx^{r-1}}\mathfrak{R}_{\delta }^{r-1}f.
\end{equation*}
\end{lemma}

Modulus of smoothness $\left\Vert \left( I-T_{\delta }\right)
^{r}f\right\Vert _{p\left( \cdot \right) }$ and \textit{K}-functional $%
K_{r}\left( f,\delta ;L^{p\left( \cdot \right) },W_{r}^{p(\cdot )}\right)
_{p\left( \cdot \right) }$ are equivalent:

\begin{theorem}
\label{reaa} If $r\in \mathbb{N}$, $p\in P^{Log}$, $f\in L^{p\left( \cdot
\right) }$, and $\delta >0$, then%
\begin{equation}
1\lesssim \frac{\left\Vert \left( I-T_{\delta }\right) ^{r}f\right\Vert
_{p\left( \cdot \right) }}{K_{r}\left( f,\delta ;L^{p\left( \cdot \right)
},W_{r}^{p(\cdot )}\right) _{p\left( \cdot \right) }}\lesssim 1  \label{real}
\end{equation}%
holds for a positive constant depend on $p,r.$
\end{theorem}

\section{Results on simultaneous approximation}

Let $\mathcal{G}_{\sigma }\left( X\right) $ be the subclass of entire
integral functions $f(z)$ of exponential type $\leq \sigma $ that belonging
to $X$ and%
\begin{equation*}
A_{\sigma }(f)_{X}:=\inf\limits_{g}\{\Vert f-g\Vert _{X}:g\in \mathcal{G}%
_{\sigma }\left( X\right) \}.
\end{equation*}%
Let $\mathcal{C}$ be the class of bounded uniformly continuous functions
defined on $\mathbb{R}$. We set $\mathcal{G}_{\sigma ,\infty }:=\mathcal{G}%
_{\sigma }\left( \mathcal{C}\right) $ and $\mathcal{G}_{\sigma ,p\left(
\cdot \right) }:=\mathcal{G}_{\sigma }\left( L^{p\left( \cdot \right)
}\right) $.

\begin{remark}
\label{rmrk}(\cite[definition given in (5.3)]{bbsv06})Let $\sigma >0$, $%
1\leq p\leq \infty $, $f\in L_{p}\left( \mathbb{R}\right) $,%
\begin{equation*}
\vartheta \left( x\right) :=\frac{2}{\pi }\frac{\sin \left( x/2\right) \sin
(3x/2)}{x^{2}}
\end{equation*}%
and%
\begin{equation*}
J\left( f,\sigma \right) =\sigma \int\nolimits_{\mathbb{R}}f\left(
x-u\right) \vartheta \left( \sigma u\right) du
\end{equation*}%
be the del\`{a} Val\`{e}e Poussin operator (\cite[definition given in (5.3)]%
{bbsv06}). It is known (see (5.4)-(5.5) of \cite{bbsv06}) that, if $f\in
L_{p}\left( \mathbb{R}\right) $, $1\leq p\leq \infty $, then,

(i) $J\left( f,\sigma \right) \in \mathcal{G}_{2\sigma }\left( L_{p}\left( 
\mathbb{R}\right) \right) $,

(ii) $J\left( g_{\sigma },\sigma \right) =g_{\sigma }$ for any $g_{\sigma
}\in \mathcal{G}_{\sigma }\left( L_{p}\left( \mathbb{R}\right) \right) $,

(iii) $\Vert J\left( f,\sigma \right) \Vert _{L_{p}\left( \mathbb{R}\right)
}\leq \frac{3}{2}\Vert f\Vert _{L_{p}\left( \mathbb{R}\right) }$,

(iv) $\left( J\left( f,\sigma \right) \right) ^{\left( r\right) }=J\left(
f^{\left( r\right) },\sigma \right) $ for any $r\in \mathbb{N}$ and $f\in
W_{p}^{r}\left( \mathbb{R}\right) $,

(v) $\Vert J\left( f,\frac{\sigma }{2}\right) -f\Vert _{L_{p}\left( \mathbb{R%
}\right) }\rightarrow 0$ (as $\sigma \rightarrow \infty $) and hence%
\begin{equation*}
\Vert \left( J\left( f,\frac{\sigma }{2}\right) \right) ^{\left( k\right)
}-f^{\left( k\right) }\Vert _{L_{p}\left( \mathbb{R}\right) }\rightarrow 0%
\text{ as }\sigma \rightarrow \infty ,
\end{equation*}%
for $f\in W_{p}^{r}\left( \mathbb{R}\right) $ and $1\leq k\leq r$,

(vi) $\left\Vert f-J\left( f,\frac{\sigma }{2}\right) \right\Vert
_{L_{p}\left( \mathbb{R}\right) }\leq \frac{5\pi }{4}\frac{4^{r}}{\sigma ^{r}%
}\Vert f^{\left( r\right) }\Vert _{L_{p}\left( \mathbb{R}\right) }$ for $%
f\in W_{p}^{r}\left( \mathbb{R}\right) .$
\end{remark}

\begin{theorem}
\label{SS1} Let $p\in P^{Log}$, $\sigma >0$, $r\in \mathbb{N}$ and $f\in
W_{r}^{p\left( \cdot \right) }$. Then%
\begin{equation}
A_{\sigma }\left( f\right) _{p\left( \cdot \right) }\lesssim \frac{1}{\sigma
^{r}}A_{\sigma }\left( f^{\left( r\right) }\right) _{p\left( \cdot \right) }
\label{obur}
\end{equation}%
holds with a positive constant depend on $p,r$.
\end{theorem}

\begin{theorem}
\label{SS2} Let $p\in P^{Log}$, $\sigma >0$, $k\in \mathbb{N}$, $r\in
\left\{ 0\right\} \cup \mathbb{N}$ and $f\in W_{r}^{p\left( \cdot \right) }$%
. Then%
\begin{equation*}
A_{\sigma }\left( f\right) _{p\left( \cdot \right) }\lesssim \Omega
_{k}\left( f,\frac{1}{\sigma }\right) _{p\left( \cdot \right) }\text{ and}
\end{equation*}%
\begin{equation}
A_{\sigma }\left( f\right) _{p\left( \cdot \right) }\lesssim \frac{1}{\sigma
^{r}}\Omega _{k}\left( f^{\left( r\right) },\frac{1}{\sigma }\right)
_{p\left( \cdot \right) }.  \label{corro}
\end{equation}
\end{theorem}

with positive constants depend on $p,k,r$.

\begin{theorem}
\label{BB3} Let $p\in P^{Log}$, $\sigma >0$ and $g_{\sigma }\in G_{\sigma
,p\left( \cdot \right) }$. Then, Bernstein's inequality 
\begin{equation*}
\Vert \left( g_{\sigma }\right) ^{\left( r\right) }\Vert _{p(\cdot
)}\lesssim \sigma ^{r}\Vert g_{\sigma }\Vert _{{p(\cdot )}}
\end{equation*}%
holds with a positive constant depend on $p,r$.
\end{theorem}

\begin{definition}
\cite[p.161]{T86}For $r$,$k\in \mathbb{N}$, $\sigma >0$, we define%
\begin{equation*}
g\left( \sigma ,r,x\right) =\left( \frac{1}{x}\sin \frac{\sigma x}{2r}%
\right) ^{2r}\text{, and}
\end{equation*}%
\begin{equation*}
G\left( \sigma ,r,k,\zeta \right) =\sum\limits_{v=1}^{k}\left( -1\right)
^{k-v}\frac{1}{v}\binom{k}{v}g\left( \sigma ,r,\frac{\zeta }{v}\right)
\end{equation*}%
For $r\geq \frac{1}{2}\left( k+2\right) $ we set%
\begin{equation*}
\gamma _{r,\sigma }:=\int\nolimits_{\mathbb{R}}\left( \frac{1}{t}\sin \frac{%
\sigma t}{2r}\right) ^{2r}dt.
\end{equation*}%
Let us introduce the Bernstein singular integral (\cite[p.161]{T86})%
\begin{equation}
D_{\sigma ,k}f(x):=\frac{\left( -1\right) ^{k+1}}{\gamma _{r,\sigma }}%
\int\nolimits_{\mathbb{R}}f(u)G\left( \sigma ,r,k,u-x\right) dt  \label{Dn}
\end{equation}%
for $r$,$k\in \mathbb{N}$, $\sigma >0$, and measurable complex valued $f$
satisfying $\int\nolimits_{\mathbb{R}}\frac{\left\vert f(u)\right\vert }{%
1+u^{2r}}du<\infty $.
\end{definition}

\begin{remark}
It is well known that, if $r$,$k\in \mathbb{N}$, $\sigma \in \left( 0,\infty
\right) $, $r\geq \frac{1}{2}\left( k+2\right) $, then $D_{\sigma ,k}f\in 
\mathcal{G}_{\sigma }\left( L^{p}\left( \mathbb{R}\right) \right) $ for $%
p\geq 1.$(\cite[p.161]{T86}).
\end{remark}

\begin{remark}
It can be shown by simple computations that%
\begin{equation*}
\gamma _{r,\sigma }=\sigma ^{2r-1}\frac{\pi ^{r}}{\left( 2r\right) ^{2r-1}}.
\end{equation*}
\end{remark}

Define $\lceil a\rceil :=\min \left\{ n\in \mathbb{N}:n\geq a\right\} $ and $%
\lfloor \sigma \rfloor :=\max \left\{ n\in \mathbb{Z}:n\leq \sigma \right\} $%
. We will take $r:=\lceil \frac{1}{2}\left( k+2\right) \rceil $ in the next
Theorem.

\begin{theorem}
\label{BernSI}Let $p\in P^{Log}$, $k\in \mathbb{N}$, $\sigma >0$, $f\in
W_{k}^{p\left( \cdot \right) }$, then%
\begin{equation}
\left\Vert f-D_{\sigma ,k}f\right\Vert _{p\left( \cdot \right) }\lesssim 
\frac{1}{\sigma ^{k}}\left\Vert f^{\left( k\right) }\right\Vert _{p\left(
\cdot \right) }  \label{tur}
\end{equation}%
holds with a positive constant depend on $p,k$.
\end{theorem}

\begin{theorem}
\label{ww}Let $p\in P^{Log}$, $k\in \mathbb{N}$, $\sigma >0$. If $f\in
L^{p\left( \cdot \right) }$, then 
\begin{equation*}
\left\Vert D_{\sigma ,k}f\right\Vert _{p\left( \cdot \right) }\lesssim
\left\Vert f\right\Vert _{p\left( \cdot \right) }
\end{equation*}%
holds with a positive constant depend on $p,k$.
\end{theorem}

\begin{corollary}
By the last Theorem \ref{ww}, if $r$,$k\in \mathbb{N}$, $\sigma \in \left(
0,\infty \right) $, $r\geq \frac{1}{2}\left( k+2\right) $, then $D_{\sigma
,k}f\in \mathcal{G}_{\sigma ,p\left( \cdot \right) }$ for $p\in P^{Log}$ and 
$f\in L^{p\left( \cdot \right) }$.
\end{corollary}

\begin{theorem}
\label{im1}Let $p\in P^{Log}$, $k\in \mathbb{N}$, $\sigma >0$. If $f\in
L^{p\left( \cdot \right) }$, then%
\begin{equation}
\left\Vert f-D_{\sigma ,k}f\right\Vert _{p\left( \cdot \right) }\lesssim
\Omega _{k}\left( f,\frac{1}{\sigma }\right) _{p\left( \cdot \right) }
\label{b3}
\end{equation}%
holds with a positive constant depending only on $k$ and $p\left( \cdot
\right) $.
\end{theorem}

\begin{theorem}
\label{S2}Let $r\in \mathbb{N}$, $p\in P^{Log}$, $\sigma >0$ and $f\in
W_{r}^{p\left( \cdot \right) }$. Then for all $k=0,1,\ldots ,r$, there
exists a positive constant depending only on $k,r$ and $p\left( \cdot
\right) $ such that%
\begin{equation*}
\left\Vert f^{\left( k\right) }-(g_{\sigma }^{\ast })^{\left( k\right)
}\right\Vert _{p\left( \cdot \right) }\lesssim \frac{1}{\sigma ^{r-k}}%
A_{\sigma }\left( f^{\left( r\right) }\right) _{p\left( \cdot \right) }
\end{equation*}%
holds for any $g_{\sigma }^{\ast }\in \mathcal{G}_{\sigma ,p\left( \cdot
\right) }$ satisfying $A_{\sigma }\left( f\right) _{p\left( \cdot \right)
}=\left\Vert f-g_{\sigma }^{\ast }\right\Vert _{p\left( \cdot \right) }$.
\end{theorem}

\begin{theorem}
\label{S3}Let $r,s\in \mathbb{N}$, $p\in P^{Log}$ and $f\in W_{r}^{p\left(
\cdot \right) }$. Then there exists a $\Phi \in \mathcal{G}_{2\sigma
,p\left( \cdot \right) }$ such that for all $k=0,1,\ldots ,r$ inequalities%
\begin{equation*}
\left\Vert f^{\left( k\right) }-\Phi ^{\left( k\right) }\right\Vert
_{p\left( \cdot \right) }\lesssim \frac{1}{\sigma ^{r-k}}\Omega _{s}\left(
f^{\left( r\right) },\frac{1}{\sigma }\right) _{p\left( \cdot \right) }
\end{equation*}%
are hold with a positive constant depending only on $k,r$ and $p\left( \cdot
\right) .$
\end{theorem}

\begin{definition}
Set $\sigma ,\eta >0,$ $f\in L^{1}\left( \mathbb{R}\right) $, $\Theta _{\eta
}f\left( x,y\right) :=f(x+\eta y)$ and%
\begin{equation*}
B_{\sigma }f(x,t):=\int\nolimits_{\mathbb{R}}\Theta _{\frac{2}{\sigma }%
}f\left( x,y\right) h\left( y,t\right) dy.
\end{equation*}
\end{definition}

\begin{remark}
The following theorem was poved in \cite{KoNa} for $\sigma =2$ with three
minor mistypes. For the sake of completeness here we will prove it when $%
\sigma >0$.
\end{remark}

\begin{theorem}
\label{Bou} Suppose that $h\left( y,t\right) ,$ $y,t\in \mathbb{R},$ is
positive measurable function with respect to y and%
\begin{equation*}
\int\nolimits_{\mathbb{R}}h\left( y,t\right) dy\lesssim 1,\quad
\int\nolimits_{\mathbb{R}}\left\vert yh_{y}^{\prime }\left( y,t\right)
\right\vert dy\lesssim 1
\end{equation*}%
with constants independent of $t$. If $\sigma >0$ and $f\in L_{1}\left( 
\mathbb{R}\right) $, then%
\begin{equation*}
\sup_{t>0}\left\vert B_{\sigma }f(\cdot ,t)\right\vert \lesssim Mf\left(
\cdot \right)
\end{equation*}%
for $t>0$ and a.e. on $\mathbb{R}$ where $Mf$ is the Hardy-Littlewood
maximal function of $f$.
\end{theorem}

\section{Proof of the results}

Let $C(A)$ be the class of continuous functions defined on $A$. For $r\in 
\mathbb{N}$, we define $C^{r}\left( A\right) $ consisting of every member $%
f\in C(A)$ such that the derivative $f^{\left( k\right) }$ exists and is
continuous on $A$ for $k=1,...,r$. We set $C^{\infty }\left( A\right)
:=\left\{ f\in C^{r}\left( A\right) \text{ for any }r\in \mathbb{N}\right\} $%
. We denote by $C_{c}\left( A\right) $, the collection of real valued
continuous functions on $A$ and support of $f$ is compact set in $A.$ We
define $C_{c}^{r}\left( A\right) :=C^{r}\left( A\right) \cap C_{c}\left(
A\right) $ for $r\in \mathbb{N}$ and $C_{c}^{\infty }\left( A\right)
:=C^{\infty }\left( A\right) \cap C_{c}\left( A\right) $. Let $L_{p}\left(
A\right) $, $1\leq p\leq \infty $ be the classical Lebesgue space of
functions on $A$.

\begin{definition}
(\cite{DHHR11}) Let $\mathbb{N}:\mathbb{=}\left\{ 1,2,3,\cdots \right\} $ be
natural numbers and $\mathbb{N}_{0}:=\mathbb{N\cup }\left\{ 0\right\} $.

(a) A family $Q$ of measurable sets $E\subset \mathbb{R}$ is called locally $%
N$-finite ($N\in \mathbb{N}$) if 
\begin{equation*}
\sum_{E\in Q}\chi _{E}\left( x\right) \leq N
\end{equation*}%
almost everywhere in $\mathbb{R}$ where $\chi _{U}$ is the characteristic
function of the set $U$.

(b) A family $Q$ of open bounded sets $U\subset \mathbb{R}$ is locally $1$%
-finite if and only if the sets $U\in Q$ are pairwise disjoint.

(c) Let $U\subset \mathbb{R}$ be \ a measurable set and%
\begin{equation*}
A_{U}f:=\frac{1}{\left\vert U\right\vert }\int\limits_{U}\left\vert f\left(
t\right) \right\vert dt.
\end{equation*}

(d) For a family $Q$ of open sets $U\subset \mathbb{R}$ we define averaging
operator by 
\begin{equation*}
T_{Q}:L_{loc}^{1}\rightarrow L^{0},
\end{equation*}%
\begin{equation*}
T_{Q}f\left( x\right) :=\sum_{U\in Q}\chi _{U}\left( x\right) A_{U}f,\quad
x\in \mathbb{R},
\end{equation*}%
where $L^{0}$ is the set of measurable functions on $\mathbb{R}$.

(e) For a measurable set $A\subset \mathbb{R}$, symbol $\left\vert
A\right\vert $ will represent the Lebesgue measure of $A$.
\end{definition}

\begin{theorem}
\label{Aver}(\cite{DHHR11})Suppose that $p\in P^{Log}$, and $f\in L^{p\left(
\cdot \right) }$. If $Q$ is $1$-finite family of open bounded subsets of $%
\mathbb{R}$ having Lebesgue measure $1$, then, the averaging operator $T_{Q}$
is uniformly bounded in $L^{p\left( \cdot \right) }$, namely,%
\begin{equation*}
\left\Vert T_{Q}f\right\Vert _{p\left( \cdot \right) }\leq \mathbf{c}%
_{4}\left\Vert f\right\Vert _{p\left( \cdot \right) }
\end{equation*}%
holds with a positive constant $\mathbf{c}_{4}$ depending only on $p$.
\end{theorem}

We define $\langle f,g\rangle :=\int\nolimits_{\mathbb{R}}f(x)g(x)dx$ when
integral exists. We will need the following Propositions.

\begin{proposition}
\label{pr1}(\cite{DHHR11})Let $p\in P^{Log}$. Then%
\begin{equation*}
\frac{1}{12\mathbf{c}_{4}}\left\Vert f\right\Vert _{p\left( \cdot \right)
}\leq \sup_{g\in L^{p^{\prime }\left( \cdot \right) }\cap C_{0}^{\infty
}:\left\Vert g\right\Vert _{p^{\prime }\left( \cdot \right) }\leq 1}\langle
\left\vert f\right\vert ,\left\vert g\right\vert \rangle \leq 2\left\Vert
f\right\Vert _{p\left( \cdot \right) }
\end{equation*}%
holds for all $f\in L^{p\left( \cdot \right) }$.
\end{proposition}

\begin{proposition}
\label{pr2}(a) $C_{c}\left( \mathbb{R}\right) $ and $C_{c}^{\infty }\left( 
\mathbb{R}\right) $ are dense subsets of $L_{p}\left( \mathbb{R}\right) $, $%
1\leq p<\infty .$(Theorems 17.10 and 23.59 of \cite[p. 415 and p. 575]{yeh}).

(b) $C_{c}\left( \mathbb{R}\right) $ contained $L_{\infty }\left( \mathbb{R}%
\right) $ but not dense (Remark 17.11 of \cite[p.416]{yeh}) in $L_{\infty
}\left( \mathbb{R}\right) .$
\end{proposition}

\begin{theorem}
\label{Fu}Let $p\in P^{Log}$. In this case,

(a) if $f\in L^{p(\cdot )}$, then, the function $F_{f}\left( \cdot \right) $
defined in (\ref{efef}) is a bounded, uniformly continuous function on $%
\mathbb{R}$,

(b) if $r\in \mathbb{N}$, and $f\in W_{r}^{p\left( \cdot \right) }$, then, $%
\left( F_{f}\right) ^{\left( k\right) }$ exists and%
\begin{equation*}
\left( F_{f}\right) ^{\left( k\right) }=F_{f^{\left( k\right) }}
\end{equation*}%
for $k\in \left\{ 1,...,r\right\} $.
\end{theorem}

\begin{proof}[\textbf{Proof of Theorem \protect\ref{Fu}}]
(a) Since $C_{0}^{\infty }$ is a dense subset of $L^{p\left( \cdot \right) }$%
, we consider functions $f\in C_{0}^{\infty }$. For any $\varepsilon >0,$
there exists $\delta :=\delta \left( \varepsilon \right) >0$ so that 
\begin{equation*}
\left\vert f\left( x+u_{1}\right) -f\left( x+u_{2}\right) \right\vert <\frac{%
\varepsilon }{1+\left\vert sptG\right\vert }
\end{equation*}%
for any $u_{1},u_{2}\in \mathbb{R}$ with $\left\vert u_{1}-u_{2}\right\vert
<\delta $. Then, there holds inequality%
\begin{equation*}
\left\vert F_{f,G}\left( u_{1}\right) -F_{f,G}\left( u_{2}\right)
\right\vert \leq \int\nolimits_{\mathbb{R}}\left\vert f\left( x+u_{1}\right)
-f\left( x+u_{2}\right) \right\vert \left\vert G\left( x\right) \right\vert
dx
\end{equation*}%
\begin{equation*}
=\int\nolimits_{sptG}\left\vert f\left( x+u_{1}\right) -f\left(
x+u_{2}\right) \right\vert \left\vert G\left( x\right) \right\vert dx
\end{equation*}%
\begin{equation*}
\leq \sup\limits_{x,u_{1},u_{2}\in sptG}\left\vert f\left( x+u_{1}\right)
-f\left( x+u_{2}\right) \right\vert \left\Vert G\right\Vert _{1,sptG}
\end{equation*}%
\begin{equation*}
\leq \frac{\varepsilon }{1+\left\vert sptG\right\vert }\left( 1+\left\vert
sptG\right\vert \right) \left\Vert G\right\Vert _{p^{\prime }\left( \cdot
\right) }\leq \varepsilon
\end{equation*}%
for any $u_{1},u_{2}\in \mathbb{R}$ with $\left\vert u_{1}-u_{2}\right\vert
<\delta $. Thus conclusion of Theorem \ref{Fu} follows. For the general case 
$f\in L_{2\pi ,\omega }^{p\left( \cdot \right) }$ there exists an $g\in
C_{0}^{\infty }$ so that%
\begin{equation*}
\left\Vert f-g\right\Vert _{p\left( \cdot \right) }<\frac{\xi }{4\left(
1+\left\vert sptG\right\vert \right) \mathbf{c}_{0}}
\end{equation*}%
for any $\xi >0$. Therefore%
\begin{equation*}
\left\vert F_{f,G}\left( u_{1}\right) -F_{f,G}\left( u_{2}\right)
\right\vert =\left\vert F_{f,G}\left( u_{1}\right) -F_{g,G}\left(
u_{1}\right) \right\vert +\left\vert F_{g,G}\left( u_{1}\right)
-F_{g,G}\left( u_{2}\right) \right\vert +
\end{equation*}%
\begin{eqnarray*}
+\left\vert F_{g,G}\left( u_{2}\right) -F_{f,G}\left( u_{2}\right)
\right\vert &=&\left\vert F_{f-g,G}\left( u_{1}\right) \right\vert +\frac{%
\xi }{2}+\left\vert F_{g-f,G}\left( u_{2}\right) \right\vert \\
&\leq &2\left( 1+\left\vert sptG\right\vert \right) \mathbf{c}_{0}\left\Vert
f-g\right\Vert _{p\left( \cdot \right) ,\omega }+\frac{\xi }{2}<\xi .
\end{eqnarray*}%
As a result $F_{f,G}$ is uniformly continuous on $\mathbb{R}$.

(b) is follow from definitions.
\end{proof}

\begin{proof}[\textbf{Proof of Theorem \protect\ref{tra}}]
Let $0\leq f,g\in L^{p(\cdot )}$. In this case there exists a constant $C>0$
such that%
\begin{align*}
\left\Vert F_{f,G}\right\Vert _{C\left( \mathbb{R}\right) }& \leq
C\left\Vert F_{g,G}\right\Vert _{C\left( \mathbb{R}\right) }=C\left\Vert
\int\nolimits_{\mathbb{R}}g\left( u+x\right) \left\vert G\left( x\right)
\right\vert dx\right\Vert _{C\left( \mathbb{R}\right) } \\
& =C\sup_{u\in \mathbb{R}}\int\nolimits_{\mathbb{R}}g\left( u+x\right)
\left\vert G\left( x\right) \right\vert dx \\
& =C\sup_{u\in sptG}\int\nolimits_{sptG}g\left( u+x\right) \left\vert
G\left( x\right) \right\vert dx \\
& \leq C\sup_{u\in sptG}\left\Vert g\left( u+\cdot \right) \right\Vert
_{1,sptG}\left\Vert G\right\Vert _{\infty }\leq C\left( 1+\left\vert
sptG\right\vert \right) \mathbf{c}_{0}\left\Vert g\right\Vert _{p(\cdot )}.
\end{align*}

On the other hand, for any $\varepsilon >0$ and appropriately chosen $\tilde{%
G}_{\varepsilon }\in L^{p^{\prime }(\cdot )}$ with 
\begin{equation*}
\int\nolimits_{\mathbb{R}}g\left( x\right) \tilde{G}_{\varepsilon }\left(
x\right) dx\geq \frac{1}{12\mathbf{c}_{4}}\left\Vert g\right\Vert _{p\left(
\cdot \right) }-\varepsilon \text{,\qquad }\left\Vert \tilde{G}_{\varepsilon
}\right\Vert _{p^{\prime }(\cdot )}\leq 1\text{,}
\end{equation*}%
(see Proposition \ref{pr1}), one can find%
\begin{eqnarray*}
\left\Vert F_{f,G}\right\Vert _{C\left( \mathbb{R}\right) } &\geq
&\left\vert F_{f,G}\left( 0\right) \right\vert \geq \int\nolimits_{\mathbb{R}%
}f\left( x\right) \left\vert G\left( x\right) \right\vert dx \\
&>&\frac{1}{12\mathbf{c}_{4}}\left\Vert f\right\Vert _{p\left( \cdot \right)
}-\varepsilon .
\end{eqnarray*}%
In the last inequality we take as $\varepsilon \rightarrow 0^{+}$ and obtain%
\begin{equation*}
\left\Vert F_{f,G}\right\Vert _{C\left( \mathbb{R}\right) }>\frac{1}{12%
\mathbf{c}_{4}}\left\Vert f\right\Vert _{p\left( \cdot \right) }\text{.}
\end{equation*}%
Combining these inequalities we get%
\begin{eqnarray*}
\left\Vert f\right\Vert _{p\left( \cdot \right) } &<&12\mathbf{c}%
_{4}\left\Vert F_{f,G}\right\Vert _{C\left( \mathbb{R}\right) }\leq 12%
\mathbf{c}_{4}C\left\Vert F_{g,G}\right\Vert _{C\left( A_{\delta }\right) }
\\
&\leq &12\mathbf{c}_{4}C\left( 1+\left\vert sptG\right\vert \right) \mathbf{c%
}_{0}\left\Vert g\right\Vert _{p\left( \cdot \right) }.
\end{eqnarray*}

For general case $f,g\in L^{p(\cdot )}$ we obtain%
\begin{equation}
\left\Vert f\right\Vert _{p\left( \cdot \right) }<24\mathbf{c}_{4}\left(
1+\left\vert sptG\right\vert \right) \mathbf{c}_{0}C\left\Vert g\right\Vert
_{p\left( \cdot \right) }.  \label{ffvv1}
\end{equation}
\end{proof}

\begin{proof}[\textbf{Proof of Theorem \protect\ref{bukun}}]
Let $0<h\leq \delta <\infty $, $p\in P^{Log}$ and $f\in L^{p\left( \cdot
\right) }$. Then, using (\ref{ffvv1}) we get%
\begin{eqnarray*}
\left\Vert \left( I-T_{h}\right) f\right\Vert _{p\left( \cdot \right) } &<&24%
\mathbf{c}_{4}\left\Vert F_{\left( I-T_{h}\right) f,G}\right\Vert _{C\left( 
\mathbb{R}\right) }\leq 24\cdot 72\mathbf{c}_{4}\left\Vert F_{\left(
I-T_{\delta }\right) f,G}\right\Vert _{C\left( \mathbb{R}\right) } \\
&\leq &1728\mathbf{c}_{4}\left( 1+\left\vert sptG\right\vert \right) \mathbf{%
c}_{0}\left\Vert \left( I-T_{\delta }\right) f\right\Vert _{p\left( \cdot
\right) }.
\end{eqnarray*}
\end{proof}

\begin{proof}[\textbf{Proof of Lemma \protect\ref{bukunA}}]
If $f\in L^{p\left( \cdot \right) }$, then, using generalized Minkowski's
integral inequality and Lemma \ref{bukun} we obtain%
\begin{equation*}
\left\Vert \left( I-\mathfrak{R}_{\delta }\right) f\right\Vert _{p\left(
\cdot \right) }=\left\Vert \frac{2}{\delta }\int\nolimits_{\delta
/2}^{\delta }\left( \frac{1}{h}\int\nolimits_{0}^{h}\left( f\left(
x+t\right) -f\left( x\right) \right) dt\right) dh\right\Vert _{p\left( \cdot
\right) }
\end{equation*}%
\begin{equation*}
=\left\Vert \frac{2}{\delta }\int\nolimits_{\delta /2}^{\delta }\left(
T_{h}f\left( x\right) -f\left( x\right) \right) dh\right\Vert _{p\left(
\cdot \right) }\leq \frac{2}{\delta }\int\nolimits_{\delta /2}^{\delta
}\left\Vert T_{\delta }f-f\right\Vert _{p\left( \cdot \right) }dh
\end{equation*}%
\begin{equation*}
\leq 1728\mathbf{c}_{4}\left( 1+\left\vert sptG\right\vert \right) \mathbf{c}%
_{0}\left\Vert T_{\delta }f-f\right\Vert _{p\left( \cdot \right) }\frac{2}{%
\delta }\int\nolimits_{\delta /2}^{\delta }dh
\end{equation*}%
\begin{equation*}
=1728\mathbf{c}_{4}\left( 1+\left\vert sptG\right\vert \right) \mathbf{c}%
_{0}\left\Vert \left( I-T_{\delta }\right) f\right\Vert _{p\left( \cdot
\right) }.
\end{equation*}
\end{proof}

\begin{proof}[\textbf{Proof of Lemma \protect\ref{lm05}}]
Using%
\begin{equation*}
\left\Vert F_{\delta \left( \mathfrak{R}_{\delta }f\right) ^{\prime
},G}\right\Vert _{C\left( \mathbb{R}\right) }=\left\Vert \delta \left(
F_{\left( \mathfrak{R}_{\delta }f\right) ,G}\right) ^{\prime }\right\Vert
_{C\left( \mathbb{R}\right) }=\delta \left\Vert \left( \mathfrak{R}_{\delta
}(F_{f,G})\right) ^{\prime }\right\Vert _{C\left( \mathbb{R}\right) }
\end{equation*}%
\begin{equation*}
\leq \cdots \leq 2\left( 37+146\ln 2^{36}\right) \left\Vert \left(
I-T_{\delta }\right) (F_{f,G})\right\Vert _{C\left( \mathbb{R}\right) }
\end{equation*}%
\begin{equation*}
=2\left( 37+146\ln 2^{36}\right) \left\Vert (F_{\left( I-T_{\delta }\right)
f,G})\right\Vert _{C\left( \mathbb{R}\right) }
\end{equation*}%
we conclude from Transference Result that%
\begin{equation*}
\delta \left\Vert (\mathfrak{R}_{\delta }f)^{\prime }\right\Vert _{p\left(
\cdot \right) }\leq \mathbf{c}_{5}\left\Vert \left( I-T_{\delta }\right)
f\right\Vert _{p\left( \cdot \right) }.
\end{equation*}%
with $\mathbf{c}_{5}:=24\mathbf{c}_{4}\left( 1+\left\vert sptG\right\vert
\right) \mathbf{c}_{0}\left( 37+146\ln 2^{36}\right) .$
\end{proof}

\begin{proof}[\textbf{Proof of Theorem \protect\ref{reaa}}]
For $r=1,2,3,\ldots $ we consider the operator%
\begin{equation*}
\mathcal{A}_{\delta }^{r}:=I-\left( I-\mathfrak{R}_{\delta }^{r}\right)
^{r}=\sum\nolimits_{j=0}^{r-1}\left( -1\right) ^{r-j+1}\binom{r}{j}\mathfrak{%
R}_{\delta }^{r\left( r-j\right) }.
\end{equation*}%
From the identity $I-\mathfrak{R}_{\delta }^{r}=\left( I-\mathfrak{R}%
_{\delta }\right) \sum\nolimits_{j=0}^{r-1}\mathfrak{R}_{\delta }^{j}$ we
find%
\begin{equation*}
\left\Vert \left( I-\mathfrak{R}_{\delta }^{r}\right) g\right\Vert _{p\left(
\cdot \right) }\leq \left( \sum\limits_{j=0}^{r-1}\mathbf{c}_{6}^{j}\right)
\left\Vert \left( I-\mathfrak{R}_{\delta }\right) g\right\Vert _{p\left(
\cdot \right) }
\end{equation*}%
with $\mathbf{c}_{6}:=24\mathbf{c}_{4}\left( 1+\left\vert sptG\right\vert
\right) \mathbf{c}_{0}.$ Therefore%
\begin{equation}
\left\Vert \left( I-\mathfrak{R}_{\delta }^{r}\right) g\right\Vert _{p\left(
\cdot \right) }\leq \left( 1728\mathbf{c}_{4}\left( 1+\left\vert
sptG\right\vert \right) \mathbf{c}_{0}\sum\limits_{j=0}^{r-1}c_{6}^{j}%
\right) \left\Vert \left( I-T_{\delta }\right) g\right\Vert _{p\left( \cdot
\right) }  \label{fdd}
\end{equation}%
\begin{equation*}
=\mathbf{c}_{7}\left\Vert \left( I-T_{\delta }\right) g\right\Vert _{p\left(
\cdot \right) }
\end{equation*}%
when $0<\delta <\infty $, $p\in P$ and $g\in L^{p\left( \cdot \right) }$.
Since $\left\Vert f-\mathcal{A}_{\delta }^{r}f\right\Vert _{p\left( \cdot
\right) }=\left\Vert \left( I-\mathfrak{R}_{\delta }^{r}\right)
^{r}f\right\Vert _{p\left( \cdot \right) }$, recursive procedure gives%
\begin{equation*}
\left\Vert f-\mathcal{A}_{\delta }^{r}f\right\Vert _{_{p\left( \cdot \right)
}}=\left\Vert \left( I-\mathfrak{R}_{\delta }^{r}\right) ^{r}f\right\Vert
_{p\left( \cdot \right) }\leq \cdots \leq \mathbf{c}_{7}^{r}\left\Vert
\left( I-T_{\delta }\right) ^{r}f\right\Vert _{p\left( \cdot \right) }.
\end{equation*}%
On the other hand, using Lemmas \ref{da} and \ref{func2}, recursively,%
\begin{equation*}
\delta ^{r}\left\Vert \frac{d^{r}}{dx^{r}}\mathfrak{R}_{\delta
}^{r}f\right\Vert _{p\left( \cdot \right) }=\delta ^{r-1}\delta \left\Vert 
\frac{d}{dx}\mathfrak{R}_{\delta }\frac{d^{r-1}}{dx^{r-1}}\mathfrak{R}%
_{\delta }^{r-1}f\right\Vert _{p\left( \cdot \right) }
\end{equation*}%
\begin{equation*}
\leq \mathbf{c}_{5}\delta ^{r-1}\left\Vert \left( I-T_{\delta }\right) \frac{%
d^{r-1}}{dx^{r-1}}\mathfrak{R}_{\delta }^{r-1}f\right\Vert _{p\left( \cdot
\right) }
\end{equation*}%
\begin{equation*}
=\mathbf{c}_{5}\delta ^{r-1}\left\Vert \left( I-T_{\delta }\right) \frac{%
d^{r-1}}{dx^{r-1}}\mathfrak{R}_{\delta }^{r-1}f\right\Vert _{p\left( \cdot
\right) }\leq \cdots \leq
\end{equation*}%
\begin{equation*}
\leq \mathbf{c}_{5}^{r-1}\delta \left\Vert \frac{d}{dx}\mathfrak{R}_{\delta
}\left( I-T_{\delta }\right) ^{r-1}f\right\Vert _{p\left( \cdot \right)
}\leq \mathbf{c}_{5}^{r}\left\Vert \left( I-T_{\delta }\right)
^{r}f\right\Vert _{p\left( \cdot \right) }.
\end{equation*}%
Thus%
\begin{equation*}
K_{r}\left( f,\delta ;L^{p\left( \cdot \right) },W_{r}^{p(\cdot )}\right)
_{p\left( \cdot \right) }\leq \left\Vert f-\mathcal{A}_{\delta
}^{r}f\right\Vert _{p\left( \cdot \right) }+\delta ^{r}\left\Vert \frac{d^{r}%
}{dx^{r}}\mathcal{A}_{\delta }^{r}f\left( x\right) \right\Vert _{p\left(
\cdot \right) }
\end{equation*}%
\begin{equation*}
\leq \mathbf{c}_{7}^{r}\left\Vert \left( I-T_{\delta }\right)
^{r}f\right\Vert _{p\left( \cdot \right) }+\sum\limits_{j=0}^{r-1}\left\vert 
\binom{r}{j}\right\vert \delta ^{r}\left\Vert \frac{d^{r}}{dx^{r}}\mathfrak{R%
}_{\delta }^{r\left( r-j\right) }f\left( x\right) \right\Vert _{p\left(
\cdot \right) }
\end{equation*}%
\begin{equation*}
\leq \mathbf{c}_{7}^{r}\left\Vert \left( I-T_{\delta }\right)
^{r}f\right\Vert _{p\left( \cdot \right) }+\mathbf{c}_{5}^{r}\sum%
\limits_{j=0}^{r-1}\left\vert \binom{r}{j}\right\vert \left\Vert \left(
I-T_{\delta }\right) ^{r}\mathfrak{R}_{\delta }^{\left( r-j\right)
}f\right\Vert _{p\left( \cdot \right) }
\end{equation*}%
\begin{equation*}
\leq \mathbf{c}_{7}^{r}\left\Vert \left( I-T_{\delta }\right)
^{r}f\right\Vert _{p\left( \cdot \right) }+\mathbf{c}_{5}^{r}\sum%
\limits_{j=0}^{r-1}\left\vert \binom{r}{j}\right\vert \mathbf{c}%
_{6}^{r-j}\left\Vert \left( I-T_{\delta }\right) ^{r}f\right\Vert _{p\left(
\cdot \right) }
\end{equation*}%
\begin{equation*}
\leq \mathbf{c}_{8}\left\Vert \left( I-T_{\delta }\right) ^{r}f\right\Vert
_{p\left( \cdot \right) }
\end{equation*}%
where%
\begin{equation*}
\mathbf{c}_{8}:=\max \left\{ \mathbf{c}_{7}^{r},\mathbf{c}%
_{5}^{r}\sum\limits_{j=0}^{r-1}\left\vert \binom{r}{j}\right\vert \mathbf{c}%
_{6}^{r-j}\right\} .
\end{equation*}

For the reverse of the last inequality, when $g\in W_{p\left( \cdot \right)
}^{r}$, we get%
\begin{equation*}
K_{r}\left( f,\delta ;L^{p\left( \cdot \right) },W_{r}^{p(\cdot )}\right)
_{p\left( \cdot \right) }\leq \left( 1+\mathbf{c}_{6}\right) ^{r}\left\Vert
f-g\right\Vert _{p\left( \cdot \right) }+\Omega _{r}\left( g,\delta \right)
_{p\left( \cdot \right) }
\end{equation*}%
\begin{equation}
\leq \left( 1+\mathbf{c}_{6}\right) ^{r}\left\Vert f-g\right\Vert _{p\left(
\cdot \right) }+2^{-r}\mathbf{c}_{6}^{r}\delta ^{r}\left\Vert g^{\left(
r\right) }\right\Vert _{p\left( \cdot \right) },  \label{mn}
\end{equation}%
and taking infimum on $g\in W_{p\left( \cdot \right) }^{r}$ in (\ref{mn}) we
obtain%
\begin{equation*}
\Omega _{r}\left( f,\delta \right) _{p\left( \cdot \right) }\leq \left( 1+%
\mathbf{c}_{6}\right) ^{r}K_{r}\left( f,\delta ;L^{p\left( \cdot \right)
},W_{r}^{p(\cdot )}\right) _{p\left( \cdot \right) }
\end{equation*}
\end{proof}

\begin{proof}[\textbf{Proof of Theorem \protect\ref{SS1}}]
The following inequality%
\begin{equation*}
A_{\sigma }\left( f\right) _{C\left( \mathbb{R}\right) }\leq \left\Vert
f-J\left( f,\frac{\sigma }{2}\right) \right\Vert _{C\left( \mathbb{R}\right)
}\leq \frac{5\pi }{4}\frac{4^{r}}{\sigma ^{r}}\Vert f^{\left( r\right)
}\Vert _{C\left( \mathbb{R}\right) }\text{,\quad }\forall f\in C^{r}(\mathbb{%
R)}
\end{equation*}%
known (see (vi) of Remark \ref{rmrk}). Now using TR we find%
\begin{equation}
\left\Vert f-J\left( f,\frac{\sigma }{2}\right) \right\Vert _{p\left( \cdot
\right) }\leq \frac{5\pi }{2}\frac{4^{r}\mathbf{c}_{6}}{\sigma ^{r}}\Vert
f^{\left( r\right) }\Vert _{p\left( \cdot \right) }\text{,\quad }\forall
f\in W_{r}^{p\left( \cdot \right) }.  \label{n}
\end{equation}

Let $r=1$. Suppose that%
\begin{equation*}
A_{\sigma }\left( f^{\prime }\right) _{p\left( \cdot \right) }=\left\Vert
f^{\prime }-g_{\sigma }^{\ast }(f^{\prime })\right\Vert _{p\left( \cdot
\right) },\ \ \ g_{\sigma }^{\ast }(f^{\prime })\in \mathcal{G}_{\sigma
,p\left( \cdot \right) }
\end{equation*}%
and%
\begin{equation*}
\digamma \left( x\right) :=\int\nolimits_{0}^{x}g_{\sigma }^{\ast
}(f^{\prime })\left( t\right) dt.
\end{equation*}%
Then $\digamma \in \mathcal{G}_{\sigma }$ (\cite[p.397]{II3}). Setting%
\begin{equation*}
\varphi \left( x\right) =f\left( x\right) -\digamma \left( x\right)
\end{equation*}%
one has%
\begin{equation*}
\left\Vert \varphi ^{\prime }\right\Vert _{p\left( \cdot \right)
}=\left\Vert f^{\prime }-g_{\sigma }^{\ast }(f^{\prime })\right\Vert
_{p\left( \cdot \right) }=A_{\sigma }\left( f^{\prime }\right) _{p\left(
\cdot \right) }.
\end{equation*}
Thus%
\begin{eqnarray*}
A_{\sigma }\left( f\right) _{p\left( \cdot \right) } &=&A_{\sigma }\left(
f-\digamma \right) _{p\left( \cdot \right) }\overset{\left( \text{\ref{n}}%
\right) }{\leq }10\pi \mathbf{c}_{6}\frac{1}{\sigma }\left\Vert \left(
f-\digamma \right) ^{\prime }\right\Vert _{p\left( \cdot \right) } \\
&=&\frac{10\pi \mathbf{c}_{6}}{\sigma }\left\Vert f^{\prime }-\digamma
^{\prime }\right\Vert _{p\left( \cdot \right) }=\frac{10\pi \mathbf{c}_{6}}{%
\sigma }\left\Vert f^{\prime }-g_{\sigma }^{\ast }(f^{\prime })\right\Vert
_{p\left( \cdot \right) } \\
&=&10\pi \mathbf{c}_{6}\frac{1}{\sigma }A_{\sigma }\left( f^{\prime }\right)
_{p\left( \cdot \right) }.
\end{eqnarray*}%
Now, result follows from the last inequality: 
\begin{eqnarray*}
A_{\sigma }\left( f\right) _{p\left( \cdot \right) } &\leq &10\pi \mathbf{c}%
_{6}\frac{1}{\sigma }A_{\sigma }\left( f^{\prime }\right) _{p\left( \cdot
\right) } \\
&\leq &\cdots \leq \left( 10\pi \mathbf{c}_{6}\right) ^{r}\frac{1}{\sigma
^{r}}A_{\sigma }\left( f^{\left( r\right) }\right) _{p\left( \cdot \right) }.
\end{eqnarray*}
\end{proof}

\begin{proof}[\textbf{Proof of Theorem \protect\ref{SS2}}]
Let $p\in P^{Log}$, $\sigma >0$, $k\in \mathbb{N}$, $r\in \left\{ 0\right\}
\cup \mathbb{N}$ and $f\in W_{r}^{p\left( \cdot \right) }$. First we
consider the case $r=0$. For every $g\in W_{k}^{p\left( \cdot \right) }$ we
find%
\begin{eqnarray*}
A_{\sigma }\left( f\right) _{p\left( \cdot \right) } &\leq &A_{\sigma
}\left( f-g\right) _{p\left( \cdot \right) }+A_{\sigma }\left( g\right)
_{p\left( \cdot \right) } \\
&\leq &\left\Vert f-g\right\Vert _{p\left( \cdot \right) }+\frac{5\pi }{2}%
\frac{4^{k}\mathbf{c}_{6}}{\sigma ^{k}}\Vert f^{\left( k\right) }\Vert
_{p\left( \cdot \right) }.
\end{eqnarray*}%
Taking infimum on $g$ in the last inequality%
\begin{equation*}
A_{\sigma }\left( f\right) _{p\left( \cdot \right) }\leq \frac{5\pi }{2}4^{k}%
\mathbf{c}_{6}K_{k}\left( f,\delta ;L^{p\left( \cdot \right)
},W_{k}^{p(\cdot )}\right) _{p\left( \cdot \right) }.
\end{equation*}%
Now using (\ref{real})%
\begin{equation*}
A_{\sigma }\left( f\right) _{p\left( \cdot \right) }\leq \mathbf{c}_{8}\frac{%
5\pi }{2}4^{k}\mathbf{c}_{6}\Omega _{k}\left( f,\frac{1}{\sigma }\right)
_{p\left( \cdot \right) }.
\end{equation*}

In the second stage we consider the case $r\in \mathbb{N}$. In this case%
\begin{equation*}
A_{\sigma }\left( f\right) _{p\left( \cdot \right) }\leq \left( 10\pi 
\mathbf{c}_{6}\right) ^{r}\frac{1}{\sigma ^{r}}A_{\sigma }\left( f^{\left(
r\right) }\right) _{p\left( \cdot \right) }
\end{equation*}%
\begin{equation*}
\leq 5\pi \mathbf{c}_{8}\left( 10\right) ^{r}\pi ^{r}\mathbf{c}%
_{6}^{r+1}2^{2k-1}\frac{1}{\sigma ^{r}}\Omega _{k}\left( f^{\left( r\right)
},\frac{1}{\sigma }\right) _{p\left( \cdot \right) }.
\end{equation*}
\end{proof}

\begin{proof}[\textbf{Proof of Theorem \protect\ref{BB3}}]
Let $p\in P^{Log}$, $\sigma >0$ and $g_{\sigma }\in \mathcal{G}_{\sigma
,p\left( \cdot \right) }$. Then, Bernstein's inequality%
\begin{equation*}
\Vert \left( g_{\sigma }\right) ^{\left( r\right) }\Vert _{C\left( \mathbb{R}%
\right) }\leq \sigma ^{r}\Vert g_{\sigma }\Vert _{C\left( \mathbb{R}\right) }%
\text{, \ \ }\forall g_{\sigma }\in \mathcal{G}_{\sigma ,\infty }
\end{equation*}%
and TR gives%
\begin{equation*}
\Vert \left( g_{\sigma }\right) ^{\left( r\right) }\Vert _{p\left( \cdot
\right) }\leq \mathbf{c}_{6}\sigma ^{r}\Vert g_{\sigma }\Vert _{p\left(
\cdot \right) }\text{, \ \ }\forall g_{\sigma }\in \mathcal{G}_{\sigma
,p\left( \cdot \right) }.
\end{equation*}
\end{proof}

\begin{proof}[\textbf{Proof of Theorem \protect\ref{BernSI}}]
Define%
\begin{equation*}
\omega _{k}\left( f,\delta \right) _{C\left( \mathbb{R}\right)
}:=\sup_{\left\vert h\right\vert \leq \delta }\left\Vert \Delta
_{t}^{k}f\right\Vert _{C\left( \mathbb{R}\right) }
\end{equation*}%
where $\Delta _{t}^{k}f\left( \cdot \right) :=\left( I-T_{h}\right)
^{r}f\left( \cdot \right) $, $T_{h}f\left( \cdot \right) :=f\left( \cdot
+h\right) $ and $I$ is the identity operator. From (\ref{Dn}), one can write%
\begin{equation*}
\left\Vert f-D_{\sigma ,k}f\right\Vert _{C\left( \mathbb{R}\right)
}=\left\Vert \frac{\left( -1\right) ^{k}}{\gamma _{r,\sigma }}\int\nolimits_{%
\mathbb{R}}\sum\limits_{v=0}^{k}\left( -1\right) ^{k-v}\binom{k}{v}f\left(
x+vt\right) g\left( \sigma ,r,t\right) dt\right\Vert _{C\left( \mathbb{R}%
\right) }
\end{equation*}%
\begin{eqnarray*}
&\leq &\frac{1}{\sigma ^{2r-1}\frac{\pi ^{r}}{\left( 2r\right) ^{2r-1}}}%
\int\nolimits_{\mathbb{R}}\left\Vert \Delta _{t}^{k}f\left( x\right)
\right\Vert _{C\left( \mathbb{R}\right) }g\left( \sigma ,r,t\right) dt \\
&\leq &\frac{\left( 2r\right) ^{2r-1}}{\pi ^{r}\sigma ^{2r-1}}\int\nolimits_{%
\mathbb{R}}\omega _{k}\left( f,t\right) _{C\left( \mathbb{R}\right) }g\left(
\sigma ,r,t\right) dt \\
&\leq &\frac{\left( 2r\right) ^{2r-1}\sigma ^{k}}{\pi ^{r}\sigma ^{2r-1}}%
\omega _{k}\left( f,\frac{1}{\sigma }\right) _{C\left( \mathbb{R}\right)
}\int\nolimits_{\mathbb{R}}\left( t+\frac{1}{\sigma }\right) ^{k}g\left(
\sigma ,r,t\right) dt \\
&\leq &\frac{\left( 2r\right) ^{2r-1}\sigma ^{k}}{\pi ^{r}\sigma ^{2r-1}}%
\frac{1}{\sigma ^{k}}\left\Vert f^{\left( k\right) }\right\Vert _{C\left( 
\mathbb{R}\right) }\int\nolimits_{\mathbb{R}}\left( t+\frac{1}{\sigma }%
\right) ^{k}g\left( \sigma ,r,t\right) dt \\
&\leq &\frac{\left( 2r\right) ^{2r-1}}{\pi ^{r}\sigma ^{2r-1}}\left\Vert
f^{\left( k\right) }\right\Vert _{C\left( \mathbb{R}\right) }\left\{ \frac{%
2^{k}}{\sigma ^{k}}\int\limits_{\left\vert t\right\vert \leq \frac{1}{\sigma 
}}\left\vert g\left( \sigma ,r,t\right) \right\vert
dt+2^{k}\int\limits_{\left\vert t\right\vert \geq \frac{1}{\sigma }%
}\left\vert t\right\vert ^{k}\left\vert g\left( \sigma ,r,t\right)
\right\vert dt\right\} .
\end{eqnarray*}%
Using $r=\lceil \frac{1}{2}\left( k+2\right) \rceil $%
\begin{equation*}
\frac{\left( 2r\right) ^{2r-1}2^{k}}{\pi ^{r}\sigma ^{2r-1}}\int_{\left\vert
t\right\vert \geq 1/\sigma }\left\vert t\right\vert ^{k}\left( \frac{1}{t}%
\sin \frac{\sigma t}{2r}\right) ^{2r}dt
\end{equation*}%
\begin{equation*}
\leq \frac{\left( 2r\right) ^{2r-1}2^{k}}{\pi ^{r}\sigma ^{2r-1}}%
\int_{\left\vert t\right\vert \geq 1/\sigma }\left( \frac{1}{t}\sin \frac{%
\sigma t}{2r}\right) ^{2r-k}dt
\end{equation*}%
\begin{equation*}
\leq \frac{\left( 2r\right) ^{2r-1}2^{k}}{\pi ^{r}\sigma ^{2r-1}}\frac{%
\sigma ^{2r-k+1}}{\left( 2r\right) ^{2r-k+1}}\int_{\mathbb{R}}\left( \frac{%
\sin u}{u}\right) ^{2}dt
\end{equation*}%
\begin{equation*}
=\frac{1}{\sigma ^{k}}\frac{2^{2k}r^{k}}{\pi ^{r}}\pi \leq \frac{1}{\sigma
^{k}}\frac{2^{2k}\left( k+2\right) ^{k}}{\pi ^{k/2}}.
\end{equation*}%
On the other hand%
\begin{equation*}
\frac{\left( 2r\right) ^{2r-1}}{\pi ^{r}\sigma ^{2r-1}}\frac{2^{k}}{\sigma
^{k}}\int_{\left\vert t\right\vert \leq 1/\sigma }\left( \frac{1}{t}\sin 
\frac{\sigma t}{2r}\right) ^{2r}dt
\end{equation*}%
\begin{equation*}
\leq \frac{\left( 2r\right) ^{2r-1}}{\pi ^{r}\sigma ^{2r-1}}\frac{2^{k}}{%
\sigma ^{k}}\int_{\mathbb{R}}\left( \frac{1}{t}\sin \frac{\sigma t}{2r}%
\right) ^{2r}dt
\end{equation*}%
\begin{equation*}
=\frac{\left( 2r\right) ^{2r-1}}{\pi ^{r}\sigma ^{2r-1}}\sigma ^{2r-1}\frac{%
\pi ^{r}}{\left( 2r\right) ^{2r-1}}=\frac{2^{k}}{\sigma ^{k}}.
\end{equation*}%
Thus%
\begin{equation*}
\left\Vert f-D_{\sigma ,k}f\right\Vert _{C\left( \mathbb{R}\right) }\leq
\left( \frac{2^{2k}\left( k+2\right) ^{k}}{\pi ^{k/2}}+2^{k}\right) \frac{1}{%
\sigma ^{k}}\left\Vert f^{\left( k\right) }\right\Vert _{C\left( \mathbb{R}%
\right) }.
\end{equation*}%
From this and TR we get%
\begin{eqnarray*}
\left\Vert f-D_{\sigma ,k}f\right\Vert _{p\left( \cdot \right) } &\leq &%
\mathbf{c}_{6}\left( \frac{2^{2k}\left( k+2\right) ^{k}}{\pi ^{k/2}}%
+2^{k}\right) \frac{1}{\sigma ^{k}}\left\Vert f^{\left( k\right)
}\right\Vert _{p\left( \cdot \right) } \\
&=&\mathbf{c}_{6}\mathbf{c}\left( k\right) \frac{1}{\sigma ^{k}}\left\Vert
f^{\left( k\right) }\right\Vert _{p\left( \cdot \right) }.
\end{eqnarray*}
\end{proof}

\begin{proof}[\textbf{Proof of Theorem \protect\ref{ww}}]
Fixed $\sigma >0,$%
\begin{equation*}
\left\Vert D_{\sigma ,k}f\right\Vert _{C\left( \mathbb{R}\right)
}=\left\Vert \frac{\left( -1\right) ^{k+1}}{\gamma _{r,\sigma }}%
\int\nolimits_{\mathbb{R}}f(u)G\left( \sigma ,r,k,u-x\right) du\right\Vert
_{C\left( \mathbb{R}\right) }
\end{equation*}%
\begin{equation*}
\left\Vert \frac{\left( -1\right) ^{k+1}}{\gamma _{r,\sigma }}\int\nolimits_{%
\mathbb{R}}\sum\limits_{v=1}^{k}\left( -1\right) ^{k-v}\binom{k}{v}f\left(
u\right) g\left( \sigma ,r,\frac{u-x}{v}\right) du\right\Vert _{C\left( 
\mathbb{R}\right) }
\end{equation*}%
\begin{equation*}
\leq \left\Vert \frac{\left( -1\right) ^{k+1}}{\gamma _{r,\sigma }}%
\int\nolimits_{\mathbb{R}}\sum\limits_{v=1}^{k}\left( -1\right) ^{k-v}\binom{%
k}{v}f\left( x+vt\right) g\left( \sigma ,r,t\right) vdt\right\Vert _{C\left( 
\mathbb{R}\right) }
\end{equation*}%
\begin{equation*}
\leq \frac{k}{\gamma _{r,\sigma }}\int\nolimits_{\mathbb{R}%
}\sum\limits_{v=1}^{k}\left\vert \binom{k}{v}\right\vert \left\Vert f\left(
x+vt\right) \right\Vert _{C\left( \mathbb{R}\right) }g\left( \sigma
,r,t\right) dt
\end{equation*}%
\begin{equation*}
\leq \left\Vert f\right\Vert _{C\left( \mathbb{R}\right)
}\sum\limits_{v=1}^{k}\left\vert \binom{k}{v}\right\vert \frac{k}{\gamma
_{r,\sigma }}\int\nolimits_{\mathbb{R}}g\left( \sigma ,r,t\right) dt\leq
k2^{k}\left\Vert f\right\Vert _{C\left( \mathbb{R}\right) }.
\end{equation*}%
Now, transference result TR gives%
\begin{equation*}
\left\Vert D_{\sigma ,k}f\right\Vert _{p\left( \cdot \right) }\leq k2^{k}%
\mathbf{c}_{6}\left\Vert f\right\Vert _{p\left( \cdot \right) }.
\end{equation*}
\end{proof}

\begin{proof}[\textbf{Proof of Theorem \protect\ref{im1}}]
We can write%
\begin{equation*}
\left\Vert f-D_{\sigma ,k}f\right\Vert _{p\left( \cdot \right) }=\left\Vert
f-\mathcal{A}_{\frac{1}{\sigma }}^{k}f+\mathcal{A}_{\frac{1}{\sigma }%
}^{k}f-D_{\sigma ,k}\mathcal{A}_{\frac{1}{\sigma }}^{k}f+D_{\sigma ,k}%
\mathcal{A}_{\frac{1}{\sigma }}^{k}f-D_{\sigma ,k}f\right\Vert _{p\left(
\cdot \right) }
\end{equation*}%
\begin{eqnarray*}
&\leq &\left\Vert f-\mathcal{A}_{\frac{1}{\sigma }}^{k}f\right\Vert
_{p\left( \cdot \right) }+\left\Vert \mathcal{A}_{\frac{1}{\sigma }%
}^{k}f-D_{\sigma ,k}\mathcal{A}_{\frac{1}{\sigma }}^{k}f\right\Vert
_{p\left( \cdot \right) }+\left\Vert D_{\sigma ,k}(\mathcal{A}_{\frac{1}{%
\sigma }}^{k}f-f)\right\Vert _{p\left( \cdot \right) } \\
&\leq &\mathbf{c}_{7}^{k}\Omega _{k}\left( f,\frac{1}{\sigma }\right)
_{p\left( \cdot \right) }+\mathbf{c}_{6}\mathbf{c}\left( k\right) \frac{1}{%
\sigma ^{k}}\left\Vert (\mathcal{A}_{\frac{1}{\sigma }}^{k}f)^{\left(
k\right) }\right\Vert _{p\left( \cdot \right) }+k2^{k}\mathbf{c}%
_{6}\left\Vert \mathcal{A}_{\frac{1}{\sigma }}^{k}f-f\right\Vert _{p\left(
\cdot \right) } \\
&\leq &\left( \mathbf{c}_{7}^{k}+\mathbf{c}_{6}\mathbf{c}\left( k\right) 
\mathbf{c}_{5}^{k}\sum\limits_{j=0}^{k-1}\left\vert \binom{k}{j}\right\vert 
\mathbf{c}_{6}^{k-j}+2^{k}k\mathbf{c}_{6}\mathbf{c}_{7}^{k}\right) \Omega
_{k}\left( f,\frac{1}{\sigma }\right) _{p\left( \cdot \right) } \\
&=&\mathbf{c}_{9}\Omega _{k}\left( f,\frac{1}{\sigma }\right) _{p\left(
\cdot \right) }
\end{eqnarray*}%
and (\ref{b3}) holds.
\end{proof}

\begin{proof}[\textbf{Proof of Theorem \protect\ref{S2}}]
Let $q\in \mathcal{G}_{\sigma }$ and $A_{\sigma }\left( f^{\left( k\right)
}\right) _{p\left( \cdot \right) }=\left\Vert f^{\left( k\right)
}-q\right\Vert _{p\left( \cdot \right) }.$ Then%
\begin{equation*}
\left\Vert f^{\left( k\right) }-(g_{\sigma }^{\ast })^{\left( k\right)
}\right\Vert _{p\left( \cdot \right) }\leq \left\Vert f^{\left( k\right)
}-\left( J\left( f,\sigma \right) \right) ^{\left( k\right) }\right\Vert
_{p\left( \cdot \right) }+\left\Vert (J\left( f,\sigma \right) )^{\left(
k\right) }-(g_{\sigma }^{\ast })^{\left( k\right) }\right\Vert _{p\left(
\cdot \right) }
\end{equation*}%
\begin{eqnarray*}
&\leq &\left\Vert f^{\left( k\right) }-q\right\Vert _{p\left( \cdot \right)
}+\left\Vert q-J\left( f^{\left( k\right) },\sigma \right) \right\Vert
_{p\left( \cdot \right) }+\left\Vert \left( J\left( f,\sigma \right)
-g_{\sigma }^{\ast }\right) ^{\left( k\right) }\right\Vert _{p\left( \cdot
\right) } \\
&\leq &A_{\sigma }\left( f^{\left( k\right) }\right) _{p\left( \cdot \right)
}+\left\Vert J\left( q-f^{\left( k\right) },\sigma \right) \right\Vert
_{p\left( \cdot \right) }+2^{k}\mathbf{c}_{6}\sigma ^{k}\left\Vert J\left(
sf,\sigma \right) -g_{\sigma }^{\ast }\right\Vert _{p\left( \cdot \right) }
\\
&\leq &\left( 1+3\mathbf{c}_{6}\right) A_{\sigma }\left( f^{\left( k\right)
}\right) _{p\left( \cdot \right) }+2^{k}\mathbf{c}_{6}\sigma ^{k}\left\Vert
J\left( f,\sigma \right) -J\left( g_{\sigma }^{\ast },\sigma \right)
\right\Vert _{p\left( \cdot \right) } \\
&\leq &\left( 1+3\mathbf{c}_{6}\right) \frac{2\mathbf{c}_{6}\left( 5\pi
4^{r-1}\right) ^{r}}{\sigma ^{r-k}}A_{\sigma }\left( f^{\left( r\right)
}\right) _{p\left( \cdot \right) }+3\mathbf{c}_{6}^{2}2^{k}\frac{2\mathbf{c}%
_{6}\left( 5\pi 4^{r-1}\right) ^{r}}{\sigma ^{r-k}}A_{\sigma }\left(
f^{\left( r\right) }\right) _{p\left( \cdot \right) } \\
&\leq &\left( 2\mathbf{c}_{6}\left( 5\pi 4^{r-1}\right) ^{r}\right) \left(
1+3\mathbf{c}_{6}+3\mathbf{c}_{6}^{2}2^{k}\right) \frac{\sigma ^{k}}{\sigma
^{r}}A_{\sigma }\left( f^{\left( r\right) }\right) _{p\left( \cdot \right) }=%
\mathbf{c}_{10}\sigma ^{k-r}A_{\sigma }\left( f^{\left( r\right) }\right)
_{p\left( \cdot \right) }
\end{eqnarray*}%
and the proof of Theorem \ref{S2} is completed.
\end{proof}

\begin{proof}[\textbf{Proof of Theorem \protect\ref{S3}}]
Let $g_{\sigma }^{\ast }\in \mathcal{G}_{\sigma }$, $A_{\sigma }\left(
f\right) _{p\left( \cdot \right) }=\left\Vert f-g_{\sigma }^{\ast
}\right\Vert _{p\left( \cdot \right) }$ and $\Phi =J\left( f,\sigma \right)
. $ Then%
\begin{equation*}
\left\Vert f-J\left( f,\sigma \right) \right\Vert _{p\left( \cdot \right)
}\leq \left\Vert f-g_{\sigma }^{\ast }+g_{\sigma }^{\ast }-J\left( f,\sigma
\right) \right\Vert _{p\left( \cdot \right) }
\end{equation*}%
\begin{equation*}
=\left\Vert f-g_{\sigma }^{\ast }+J\left( g_{\sigma }^{\ast },\sigma \right)
-J\left( f,\sigma \right) \right\Vert _{p\left( \cdot \right) }
\end{equation*}%
\begin{equation*}
\leq A_{\sigma }\left( f\right) _{p\left( \cdot \right) }+3\mathbf{c}%
_{6}\left\Vert f-g_{\sigma }^{\ast }\right\Vert _{p\left( \cdot \right)
}=\left( 1+3\mathbf{c}_{6}\right) A_{\sigma }\left( f\right) _{p\left( \cdot
\right) }
\end{equation*}%
and%
\begin{equation*}
\left\Vert f-J\left( f,\sigma \right) \right\Vert _{p\left( \cdot \right)
}\leq \left( 1+3\mathbf{c}_{6}\right) A_{\sigma }\left( f\right) _{p\left(
\cdot \right) }
\end{equation*}%
\begin{equation*}
\leq \left( 1+3\mathbf{c}_{6}\right) 5\pi \mathbf{c}_{8}\left( 10\right)
^{r}\pi ^{r}\mathbf{c}_{6}^{r+1}2^{2s-1}\frac{1}{\sigma ^{r}}\Omega
_{s}\left( f^{\left( r\right) },1/\sigma \right) _{p\left( \cdot \right) }.
\end{equation*}%
Now, from%
\begin{equation*}
\left\Vert f-g_{\sigma }^{\ast }\right\Vert _{p\left( \cdot \right) }\leq 
\frac{\pi \mathbf{c}_{8}\left( 10\right) ^{r}\pi ^{r}\mathbf{c}%
_{6}^{r+1}2^{2s-1}}{\sigma ^{r}}\Omega _{s}\left( f^{\left( r\right) },\frac{%
1}{\sigma }\right) _{p\left( \cdot \right) }
\end{equation*}%
we obtain%
\begin{equation*}
\left\Vert J\left( f,\sigma \right) -g_{\sigma }^{\ast }\right\Vert
_{p\left( \cdot \right) }\leq \frac{\mathbf{c}_{11}}{\sigma ^{r}}\Omega
_{s}\left( f^{\left( r\right) },\frac{1}{\sigma }\right) _{p\left( \cdot
\right) }
\end{equation*}%
with%
\begin{equation*}
\mathbf{c}_{11}=\pi \mathbf{c}_{8}\left( 10\right) ^{r}\pi ^{r}\mathbf{c}%
_{6}^{r+1}2^{2s-1}\left( \left( 1+3\mathbf{c}_{6}\right) 5+1\right) .
\end{equation*}%
Hence 
\begin{eqnarray*}
\left\Vert f^{\left( k\right) }-\left( J\left( f,\sigma \right) \right)
^{\left( k\right) }\right\Vert _{p\left( \cdot \right) } &\leq &\left\Vert
f^{\left( k\right) }-\left( g_{\sigma }^{\ast }\right) ^{\left( k\right)
}\right\Vert _{p\left( \cdot \right) }+\left\Vert \left( J\left( f,\sigma
\right) \right) ^{\left( k\right) }-\left( g_{\sigma }^{\ast }\right)
^{\left( k\right) }\right\Vert _{p\left( \cdot \right) } \\
&\leq &\mathbf{c}_{10}\sigma ^{k-r}A_{\sigma }\left( f^{\left( r\right)
}\right) _{p\left( \cdot \right) }+2^{k}\mathbf{c}_{6}\sigma ^{k}\frac{%
\mathbf{c}_{11}}{\sigma ^{r}}\Omega _{s}\left( f^{\left( r\right) },\frac{1}{%
\sigma }\right) _{p\left( \cdot \right) } \\
&\leq &\left( \mathbf{c}_{10}\frac{5\pi \mathbf{c}_{8}}{2}4^{s}\mathbf{c}%
_{6}+2^{k}\mathbf{c}_{6}\mathbf{c}_{11}\right) \sigma ^{k-r}\Omega
_{s}\left( f^{\left( r\right) },1/\sigma \right) _{p\left( \cdot \right) }
\end{eqnarray*}%
and the proof is completed.
\end{proof}

\begin{proof}[\textbf{Proof of Theorem \protect\ref{Bou}}]
Given $x\in \mathbb{R}$, let%
\begin{equation*}
\Gamma \left( y\right) :=\int\nolimits_{0}^{y}\Theta _{\frac{2}{\sigma }%
}f\left( x,u\right) du,\quad y>0,
\end{equation*}%
and $a,b>0$. Integration by parts gives%
\begin{equation*}
\int\nolimits_{-a}^{b}\Theta _{\frac{2}{\sigma }}f\left( x,y\right) h\left(
y,t\right) dy=\int\nolimits_{-a}^{b}h\left( y,t\right) d\Gamma \left(
y\right)
\end{equation*}%
\begin{equation*}
=\Gamma \left( y\right) h\left( y,t\right) \mid
_{-a}^{b}-\int\nolimits_{-a}^{b}h_{y}^{\prime }\left( y,t\right) \Gamma
\left( y\right) dy.
\end{equation*}%
Since $\Gamma \left( y\right) \leq \left\vert y\right\vert Mf\left( x\right) 
$ we obtain%
\begin{equation*}
\left\vert \int\nolimits_{-a}^{b}\Theta _{\frac{2}{\sigma }}f\left(
x,y\right) h\left( y,t\right) dy\right\vert \leq Mf\left( x\right) \left(
\int\nolimits_{-a}^{b}\left\vert yh_{y}^{\prime }\left( y,t\right)
\right\vert dy+h\left( y,t\right) \mid _{-a}^{b}\right) .
\end{equation*}%
Now%
\begin{equation*}
\mathbf{c}_{12}\geq \int\nolimits_{\mathbb{R}}h\left( y,t\right) dy\geq
\int\nolimits_{-a}^{b}h\left( y,t\right) dy=h\left( y,t\right) \mid
_{-a}^{b}-\int\nolimits_{-a}^{b}yh_{y}^{\prime }\left( y,t\right) dy
\end{equation*}%
gives%
\begin{equation*}
\left\vert \int\nolimits_{-a}^{b}\Theta _{\frac{2}{\sigma }}f\left(
x,y\right) h\left( y,t\right) dy\right\vert \leq \left( \mathbf{c}_{12}+2%
\mathbf{c}_{13}\right) Mf\left( x\right)
\end{equation*}%
for any $t>0.$ The last inequality implies the result.
\end{proof}

\end{document}